\documentclass[a4paper,10pt]{article}

\usepackage[english]{babel}
\usepackage[utf8]{inputenc}

\usepackage{amsmath}
\usepackage{amsfonts}
\usepackage{amssymb}
\usepackage{amsthm}
\usepackage{graphicx}
\usepackage{pdfpages}

\newcommand{\mytext}[1]{ \: \textrm{#1} \: }
\newcommand{\mysetdescr}[2]{\left\{ #1 \: \left| \: #2 \right. \right\} }

\newcommand{\inverse}[1]{{#1}^{-1}}

\def\G{{\cal G}}
\def\H{{\cal H}}
\def\L{{\cal L}}
\def\P{{\cal P}}

\newcommand{\mf}[1]{\mathfrak{ #1 }}
\newcommand{\fC}{\mf{C}}
\newcommand{\fF}{\mf{F}}
\newcommand{\fG}{\mf{G}}
\newcommand{\fI}{\mf{I}}
\newcommand{\fM}{\mf{M}}
\newcommand{\fN}{\mf{N}}
\newcommand{\fP}{\mf{P}}
\newcommand{\fU}{\mf{U}}

\newcommand{\Zy}{{Z \cup \{ y \}}}

\newcommand{\UXY}{\fU ( X, Y )}
\newcommand{\CQXY}{\fC_Q ( X, Y )}

\newcommand{\MQXY}{\fM_Q ( X, Y )}
\newcommand{\MstQXY}{\fM^*_Q ( X, Y )}

\newcommand{\MQdXY}{\fM_{Q^d} ( X, Y )}

\newcommand{\FQXY}{\fF_Q ( X, Y )}
\newcommand{\FstQXY}{\fF^*_Q ( X, Y )}

\newcommand{\MstQXZy}{\fM^*_{Q + A_y} \left( X, \Zy \right)}

\newcommand{\CQdXZ}{\fC_{Q^d} ( X, Z )}

\newcommand{\FstQXZy}{\fF^*_{Q + A_y} \left( X, \Zy \right)}

\newcommand{\setyps}{\{ y \}}

\def\BP{\begin{proof}}
\def\EP{\end{proof}}

\newcommand{\equivA}{ \; \equiv \; }

\newcommand{\eqA}{ \; = \; }
\newcommand{\inA}{ \in \; }
\newcommand{\subseteqA}{ \subseteq \; }
\newcommand{\supseteqA}{ \supseteq \; }
\newcommand{\capA}{ \cap \; }
\newcommand{\cupA}{ \cup \; }

\begin{document}

\theoremstyle{plain}
\newtheorem{condition}{Condition}
\newtheorem{theorem}{Theorem}
\newtheorem{definition}{Definition}
\newtheorem{corollary}{Corollary}
\newtheorem{lemma}{Lemma}
\newtheorem{proposition}{Proposition}

\title{\bf A generalization of a theorem of Ern\'{e}}
\author{\sc Frank a Campo}
\date{\small Viersen, Germany\\
{\sf acampo.frank@gmail.com} \\
May 2021}

\maketitle

\begin{abstract}
\noindent Let $X$ be a finite set, $Z \subseteq X$ and $y \notin X$. Marcel Ern\'{e} showed in 1981, that the number of posets on $X$ containing $Z$ as an antichain equals the number of posets $R$ on $X \cup \setyps$ in which the points of $\Zy$ are exactly the maximal points of $R$. We prove the following generalization: For every poset $Q$ with carrier $Z$, the number of posets on $X$ containing $Q$ as an induced sub-poset equals the number of posets $R$ on $X \cup \setyps$ which contain $Q^d + A_y$ as an induced sub-poset and in which the maximal points of $Q^d + A_y$ are exactly the maximal points of $R$. Here, $Q^d$ is the dual of $Q$, $A_y$ is the singleton-poset on $y$, and $Q^d + A_y$ denotes the direct sum of $Q^d$ and $A_y$.
\newline

\noindent{\bf Mathematics Subject Classification:}\\
Primary: 06A07. Secondary: 06A06.\\[2mm]
{\bf Key words:} poset, enumeration of posets.
\end{abstract}

\section{Introduction} \label{sec_Introduction}

The determination of the number of partially ordered sets (posets) on $n$ points is one of the problems in order theory still far from being solved. The number equals the number of $T_0$-topologies on $n$ points, because according to the well-known results of Alexandroff \cite{Alexandroff_1935,Alexandroff_1937} and Birkhoff \cite{Birkhoff_1967}, these structures are in one-to-one correspondence. The task to count posets on $n$ points can in a natural way be broken down into the sub-tasks of counting posets providing additional properties, e.g., a fixed number $k$ of down-sets (or antichains), or - equivalently - a fixed number $k$ of open (or closed) sets of $T_0$-topologies on $n$ points. Early attempts are Sharp \cite{Sharp_1966,Sharp_1968}, Stanley \cite{Stanley_1974},  Butler and Markowsky \cite{Butler_1972,Butler_Markowsky_1973}, Ern\'{e} \cite{Erne_1972,Erne_1974}, and Parchmann \cite{Parchmann_1975}. In 1971, Stanley \cite{Stanley_1971} determined all non-isomorphic posets with a down-set-number $k \geq \frac{7}{16} \cdot 2^n$.

Ern\'{e} \cite{Erne_1981} presented in 1981 combinatorial identities and showed that the number of posets on $n$ points with a fixed antichain of $m$ elements equals the number of posets on $n+1$ points with fixed $m+1$ maximal elements - a result we will generalize in this article. Vollert \cite[Satz 2.18]{Vollert_1987} proved in 1987 a striking regularity in the binary representation of the relative number of down-sets of posets.

In \cite{Erne_Stege_1991}, Ern\'{e} and Stege improved an algorithm for counting posets and determined their number for $n \leq 14$; this article contains also a survey over the development of the art of counting up to 1991. In \cite{Erne_1992}, Ern\'{e} published results about the number of posets with a certain structure (see also Culberson and Rawlins \cite{Culberson_Rawlins_1991}), and in \cite{Erne_Stege_1995}, Ern\'{e} and Stege presented a general combinatorial recursive approach and applied it to posets. In 2000, Heitzig and Reinhold \cite{Heitzig_Reinhold_2000} calculated with an improved algorithm the numbers of non-isomorphic ("unlabelled") posets for $n = 13$ and the number of all posets including isomorphic ("labelled") ones for $n = 15$; in \cite{Heitzig_Reinhold_2002}, they counted lattices for $n \leq 18$. In 2002, Ern\'{e}, Heitzig, and Reinhold \cite{Erne_etal_2002} investigated the numbers of distributive lattices with $n$ elements and of posets with $n$ antichains for $n < 50$.

In 2002, Brinkmann and McKay \cite{Brinkmann_McKay_2002} developed a new algorithm and applied it to the numbers of unlabelled posets and $T_0$-topologies for $n = 15, 16$ and of labelled posets for $n = 18, 19$. Pfeiffer \cite{Pfeiffer_2004} calculated in 2004 the numbers of several transitive relations for $n \leq 12$. Brinkmann and McKay \cite{Brinkmann_McKay_2005} extended Pfeiffer's approach to $n = 15, 16$.

Dershowitz and Zaks \cite{Dershowitz_Zaks_1980} enumerated ordered trees in 1980, and Bayoumi, El-Zahar, and Khamis \cite{Bayoumi_etal_1994} counted two-dimensional posets in 1994. In 2006-2014, Benoumhani and Kolli \cite{Benoumhani_2006,Kolli_2007,Kolli_2014,Benoumhani_Kolli_2010} determined the numbers of labelled and unlabelled $T_0$-topologies with $n$ points and $k$ open sets for many combinations of $n$ and $k$ by careful case-by-case studies. Kolli presented in  \cite{Kolli_2007} all $T_0$-topologies for $k \geq \frac{3}{8} \cdot 2^n$ and in \cite{Kolli_2014} for $k \geq \frac{5}{16} \cdot 2^n$. By means of a new systematic approach, the author \cite{aCampo_2018} was able to present all non-isomorphic posets with down-set-number $k \geq \frac{1}{4} \cdot 2^n$ in 2018, including the posets Vollert \cite{Vollert_1987} had pointed at. Ragnarsson and Tenner \cite{Ragnarsson_Tenner_2010} studied in 2010 the smallest possible number $n$ of points of a topology having $k$ open sets. In 2019, the author and Ern\'{e} \cite{aCampo_Erne_2019} presented exponential functions for the number of extensions of finite posets with a fixed set of minimal points.

The subject of the present article is a generalization of the difficult result proven by Marcel Ern\'{e} \cite{Erne_1981} in 1981:

\begin{theorem}[\cite{Erne_1981} Theorem 3.4] \label{theo_Erne}
Let $X$ and $Z$ be disjoint finite sets and let $y$ be a point not contained in $X \cup Z$. The following numbers are equal:
\begin{itemize}
\item The number of posets on $X \cup Z$ with $Z$ being an antichain;
\item the number of posets $R$ on $X \cup \Zy$ in which the points of $\Zy$ are exactly the maximal points of $R$.
\end{itemize}
\end{theorem}
We generalize this theorem as follows:
\begin{theorem} \label{theo_VVErne_Intro}
Let $X$ and $Z$ be disjoint finite sets, let $Q$ be a poset on $Z$, and let $y$ be a point not contained in $X \cup Z$. The following numbers are equal:
\begin{itemize}
\item The number of posets on $X \cup Z$ containing $Q$ as an induced sub-poset;
\item the number of posets $R$ on $X \cup \Zy$ for which $Q^d + A_y$ is an induced sub-poset and in which the maximal points of $Q^d + A_y$ are exactly the maximal points of $R$.
\end{itemize}
Here, $Q^d$ is the dual of $Q$, $A_y$ is the singleton-poset on $y$, and $Q^d + A_y$ denotes the direct sum of $Q^d$ and $A_y$.
\end{theorem}

After recalling common terms in Section \ref{sec_notation}, we develop the proof of Theorem \ref{theo_VVErne_Intro} in two steps. The first one is done in Section \ref{sec_Verne} where we reach in Theorem \ref{theo_Verne} a result corresponding to Theorem \ref{theo_VVErne_Intro}, but still with the restriction that the set $Z$ has to be {\em convex} in $R$. In Section \ref{sec_VVerne}, we skip this restriction and finish with Theorem \ref{theo_VVErne} the proof of Theorem \ref{theo_VVErne_Intro}. In the Figures \ref{figure_Beispiel_Az}-\ref{figure_Beispiel_V} at the end of the article, Theorem \ref{theo_VVErne} is illustrated.

In the proof of Theorem \ref{theo_Verne} in Section \ref{sec_Verne}, we follow in detail the ingenious original proof of Ern\'{e}'s Theorem in \cite{Erne_1981}, just making the required modifications at every point. We decided to present this demanding etude in full length, because (in our eyes) every substantial reduction by ``just as in \cite{Erne_1981}, we see ...'' would open a gap difficult to bridge for the reader. Section \ref{sec_VVerne} has no counterpart in \cite{Erne_1981}.

\section{Basics and Notation} \label{sec_notation}

Given a set $X$, a partial order relation $R$ on $X$ is a reflexive, antisymmetric, and transitive subset of $X \times X$. The pair $P = (X,R)$ is called a {\em partially ordered set} or simply a {\em poset}; $X$ is called the {\em carrier} of $P$. Because every partial order relation on $X$ is reflexive, the {\em diagonal} $\Delta_X \equiv \mysetdescr{(x,x)}{x\in X}$ is always contained in it. Given the partial order relation $R$ of a poset $P$, we can thus reconstruct the carrier by means of $R$; we denote it by $c(R)$: $P = (c(R),R)$.

Due to this one-to-one correspondence between the posets with carrier $X$ and the partial order relations on $X$, there is no difference between dealing with posets and dealing with partial order relations. And indeed, the coming sections work exclusively with partial order relations. Therefore, we direct all definitions and terms towards partial order relations instead of posets. The symbol $\fP(X)$ denotes the set of all partial order relations with carrier $X$, and the abbreviation ``p.o.r.'' means ``partial order relation''.

On any set $X$, we can define two elementary partial order relations: The {\em antichain} $\Delta_X$ and the {\em chain} $C_X$ which is up to isomorphism characterized by $(x,y) \in C_X$ or $(y,x) \in C_X$ for all $x, y \in X$.

Given a p.o.r.\ $R$ with carrier $X$ and a subset $M \subseteq X$, the {\em induced p.o.r.} $R \vert_M$ is defined as $R \cap (M \times M)$; the carrier of $R \vert_M$ is $M$. Given partial order relations $R \in \fP(X)$ and $Q \in \fP(Y)$, we call $Q$ a sub-p.o.r. of $R$ iff $Q \subseteq R$ (which implies $Y \subseteq X$). An induced p.o.r.\ is always a sub-p.o.r.

For a partial order relation $R$ on $X$, the {\em dual} $R^d$ is defined by 
\begin{align*}
R^d & \equivA \mysetdescr{ (y,x) \in X \times X }{ (x,y) \in R }.
\end{align*}
The diagram of the poset corresponding to $R^d$ is the diagram of the poset corresponding to $R$ turned upside-down.

For two partial order relations $R_1$ on $X_1$ and $R_2$ on $X_2$ with $X_1 \cap X_2 = \emptyset$, their {\em direct sum} is the p.o.r.\ $R_1 + R_2 \equiv R_1 \cup R_2$. We get the diagram of the corresponding poset by juxtaposing the diagrams of the posets $P_1 \equiv (X_1,R_1)$ and $P_2 \equiv (X_2,R_2)$ without any connection between them. The {\em ordinal sum} of $R_1$ and $R_2$ is the p.o.r.\ $R_1 \oplus R_2 \equiv R_1 \cup R_2 \cup ( X_1 \times X_2 )$; here, we get the corresponding diagram by putting $P_2$ on top of $P_1$ and connecting every maximal point of $P_1$ with every minimal point of $P_2$.

Let $R \in \fP(X)$. A subset $L \subseteq X$ is called a {\em lower end} of $R$ iff we have $x \in L$ for every $x \in X$ with $(x,y) \in R$ for a $y \in L$. The set of lower ends of $R$ is denoted by $\L(R)$. For every subset $M \subseteq X$ and every $x \in X$, we define the lower ends created by $M$ and $x$, respectively, as
\begin{align*}
\downarrow_R M & \equivA \mysetdescr{ x \in X }{ (x,y) \in R \mytext{ for a } y \in M }, \\
\downarrow_R x & \equivA \downarrow_R \{ x \}.
\end{align*}
The inclusion ``$\subseteq$'' introduces a partial order relation on the set $\L(R)$ making it a lattice. Due to the representation theorem of Birkhoff \cite[p.\ 59]{Birkhoff_1967}, a p.o.r.\ $R \in \fP(X)$ with finite carrier $X$ can uniquely be identified by means of $\L(R)$.

{\em Upper ends} are the dual of lower ends: A subset $U \subseteq X$ is called an {\em upper end} of $R$ iff we have $x \in U$ for every $x \in X$ with $(y,x) \in R$ for a $y \in U$. Correspondingly, we define for every subset $M \subseteq X$ and every $x \in X$,
\begin{align*}
\uparrow_R M & \equivA \mysetdescr{ x \in X }{ (y,x) \in R \mytext{ for a } y \in M }, \\
\uparrow_R x & \equivA \uparrow_R \{ x \}.
\end{align*}
A subset $M$ of $X$ is a lower end of $R$, iff $R \cap ( ( X \setminus M ) \times M ) = \emptyset$; correspondingly, it is an upper end iff $R \cap ( M \times ( X \setminus M ) ) = \emptyset$.

Given $R \in \fP(X)$, a point $x \in X$ is called a {\em maximal point} of $R$, iff $\uparrow_R x = \{ x \}$. $\max R$ denotes the set of maximal points of $R$. {\em Minimal points} of $R$ are defined correspondingly. For a subset $M \subseteq X$, we define $\max M \equiv \max R \vert_M$.

Given partial order relations $P \in \fP(X)$ and $Q \in \fP(Y)$, $X$ and $Y$ not necessary disjoint, the {\em transitive hull} of $P$ and $Q$ is the smallest transitive relation $T$ on $X \cup Y$ with $R \cup Q \subseteq T$. If $P$ and $Q$ are both reflexive, then also $T$ is, but antisymmetry of $P$ and $Q$ is not necessary inherited by $T$. However, in the case discussed after Definition \ref{def_GQAY} in Section \ref{sec_VVerne}, $T$ is antisymmetric, too, hence a partial order relation.

We call a subset $M \subseteq X$ {\em convex} in $R$, iff $(a,x), (x,b) \in R$ implies $x \in M$ for all $a, b \in M$ and all $x \in X$. A lower end or upper end of $R$ is always convex.

The {\em convex hull} of a set $M \subseteq X$ is the smallest convex set in $R$ containing $M$. In Section \ref{sec_VVerne}, we need:
\begin{lemma} \label{lemma_ConstrConv}
Let $X$ be any finite set, $R \in \fP( X )$, and $M \subseteq X$. With
\begin{align*}
\gamma_R(M) & \equivA \bigcup_{(m,n) \in M \times M} ( \uparrow_R m ) \cap ( \downarrow_R n ),
\end{align*}
$\gamma_R(M)$ is the convex hull of $M$ in $R$, and we have $\max M = \max \gamma_R(M)$. Additionally, if $\max M = \max R$, then $\gamma_R(M)$ is an upper end of $R$.
\end{lemma}
\BP Let $a, a' \in \gamma_R(M)$ and $x \in X$ with $(a,x), (x,a') \in R$. There exist $m, n, m', n' \in M$ with $(m,a), (a,n), (m',a'), (a',n') \in R$, hence $(m,x), (x,n') \in R$, and $x \in \gamma_R(M)$ follows. $\gamma_R(M)$ is thus a convex subset of $R$ with $M \subseteq \gamma_R(M)$. Conversely, every convex set $C$ in $R$ containing $M$ must contain all the sets $( \uparrow_R m ) \cap ( \downarrow_R n )$, $m, n \in M$, thus $\gamma_R(M) \subseteq C$.

Let $x \in \max \gamma_R(M)$. Because of $x \in \gamma_R(M)$, there exist $m, n \in M$ with $(m,x), (x,n) \in R$, thus $x = n \in M$. But then $x \in \max M$ due to $M \subseteq \gamma_R(M)$. Conversely, let $m \in \max M$. Because $\gamma_R(M)$ is finite, there exists a $x \in \max \gamma_R(M)$ with $(m,x) \in R$. We have already seen $x \in M$, thus $m = x \in \max \gamma_R(M)$.

Assume $\max R = \max M$, and let $a \in \gamma_R(M)$ and $x \in X$ with $(a,x) \in R$. There exist $m, n \in M$ with $(m,a), (a,n) \in R$, hence $(m,x) \in R$. In the case of $x \notin \gamma_R(M)$, there exists no $n' \in M$ with $(x,n') \in R$, hence $\uparrow_R x \subseteq X \setminus M$ in contradiction to $\max R = \max M$.

\EP

In order to avoid unnecessary formalism, we apply the terms ``lower end'', ``upper end'', ``convex'' etc.\ defined for a subset $M$ of $X$ also to the partial order relation induced by $M$. We call thus an induced partial order relation $Q = R \vert_M$ a lower end, an upper end, or convex in $R$, iff the carrier $M$ of $Q$ has the respective property in $R$.

\section{A first generalization of Ern\'{e}'s Theorem} \label{sec_Verne}

In this section, we reach in Theorem \ref{theo_Verne} a first generalization of Ern\'{e}'s Theorem, focusing on partial order relations $R$ containing $Q$ as an induced {\em convex} sub-p.o.r. As already announced in the introduction, we follow in detail the line of thought in the original proof Theorem \ref{theo_Erne} in \cite{Erne_1981}. Convexity of $Q$ in $R$ is needed in the proofs of Corollary \ref{coro_UXY_C} and Lemma \ref{lemma_tau}, which in turn are used in the proof of Theorem \ref{theo_Verne} in showing that a certain mapping is bijective. The fact that $Z$ is not anymore an antichain affects directly the proofs of Lemma \ref{lemma_tau} (needed for Corollary \ref{coro_MstQXY}) and of Theorem \ref{theo_Verne}.

\begin{definition} \label{def_MQXY}
For finite disjoint sets $X$ and $Y$ and every $Q \in \fP(Y)$,
\begin{align*}
\UXY & \equivA \mysetdescr{ R \in \fP(X \cup Y) }{ Y \mytext{is an upper end of } R }, \\
\CQXY & \equivA \mysetdescr{ R \in \fP(X \cup Y) }{ R \vert_Y = Q \mytext{ is convex in } R}, \\
\MQXY & \equivA \UXY \cap \CQXY, \\
\MstQXY & \equivA \mysetdescr{ R \in \CQXY }{ \max Q = \max R }.
\end{align*}
Additionally, we define
\begin{align*}
\FQXY & \equivA \bigcup_{P \in \fP(X)} \H_Q( Y, \L(P) ),
\end{align*}
where $\H_Q(Y, \L(P))$ is the set of mappings $f$ from $Y$ to $\L(P)$ with $(a,b) \in Q \Rightarrow f(a) \subseteq f(b)$ for all $a,b \in Y$.
\end{definition}
For every $R \in \MQXY$, we have $\max Q \subseteq \max R$, and due to Lemma \ref{lemma_ConstrConv}, $\MstQXY$ is a subset of $\MQXY$. The elements of $\FQXY$ are triplets $f = (Y,F,\L(P))$ with certain relations $F \subseteq Y \times \L(P)$ and $P \in \fP(X)$. Because every partial order relation with finite carrier is uniquely determined by its lower end lattice, we can identify $P$ by means of the third component of $f$, and we define
\begin{align*}
S(f) & \equivA P.
\end{align*}
\begin{lemma} \label{lemma_MstQXY}
The mapping
\begin{align} \label{def_Phi}
\begin{split}
\Phi : \FQXY & \rightarrow \MQXY \\
f & \mapsto S(f) \; \cupA Q \; \cupA \bigcup_{y \in Y} ( f(y) \times \{ y \} )
\end{split}
\end{align}
is a bijection with inverse
\begin{align} \label{inv_Phi}
\begin{split}
\inverse{\Phi}(R) : Y & \rightarrow \L( R \vert_X ), \\
y & \mapsto \left( \downarrow_R y \right) \setminus Y.
\end{split}
\end{align}
\end{lemma}
\BP Let $f \in \FQXY$.  Because of $\Phi(f) \cap ( Y \times X ) = \emptyset$, we have $y \in X \Rightarrow x \in X$ for all $(x,y) \in \Phi(f)$.

Obviously, $\Phi( f )$ is reflexive. Let $(x,y), (y,x) \in \Phi(f)$. If $y \in X$ or $x \in X$, then $x,y \in X$, and $x = y$ follows. And in the remaining case $x,y \in Y$, we have $(x,y), (y,x) \in Q$, hence $x = y$.

Let $(a,b), (b,c) \in \Phi(f)$. There are four possible cases:
\begin{itemize}
\item $c \in X$: Then $(a,b), (b,c) \in S(f)$, hence $(a,c) \in S(f) \subseteq \Phi( f )$.
\item $c \in Y$, $b \in X$: Then $(a,b) \in S(f)$ and $b \in f(c) \in \L(S(f))$, hence $a \in f(c)$, and $(a,c) \in \Phi(f)$ follows.
\item $b, c \in Y$, $a \in X$: Then $a \in f(b) \in \L(S(f))$. Due to $f \in \H_Q(Y, \L(S(f)))$, we have $f(b) \subseteq f(c)$, and $(a,c) \in \Phi(f)$ follows.
\item $a, b, c \in Y$: Then $(a,b), (b,c) \in Q$, thus $(a,c) \in Q \subseteq S(f)$.
\end{itemize}
Therefore, $\Phi(f)$ is a partial order relation on $X$. Because of $\Phi(f) \cap ( Y \times X ) = \emptyset$, the set $Y$ is an upper end in $\Phi(f)$, thus a convex subset with $\Phi(f) \vert_Y = Q$, and $\Phi(f) \in \MQXY$ is shown.

Let $R \in \MQXY$. Because $Y$ is an upper end in $R$, the sets $\left( \downarrow_R y \right) \setminus Y$ with $y \in Y$ are lower ends in $R \vert_X$. The mapping $\psi$ described in \eqref{inv_Phi} is thus an element of $\H_Q(Y,\L(R \vert_X )) \subseteq \FQXY$, and due to $\left( \left( \downarrow_R y \right) \cap Y \right) \times \setyps \subseteq R \vert_Y = Q$ for all $y \in Y$, we have 
\begin{align*}
\Phi(\psi) & \eqA R \vert_X \; \cupA R \vert_Y \; \cupA ( R \cap ( ( X \cup Y ) \times Y )),
\end{align*}
and $\Phi( \psi ) = R$ follows because $Y$ is an upper end of $R$. Therefore, $\Phi$ is bijective.

\EP

\begin{corollary} \label{coro_MstQXY}
With
\begin{align*}
\FstQXY & \equivA \mysetdescr{ f \in \FQXY }{ X = \cup f[ Z ] },
\end{align*}
the mapping $\Phi$ induces a bijection between $\FstQXY$ and $\MstQXY$.
\end{corollary}
\BP Let $f \in \FstQXY$. According to Lemma \ref{lemma_MstQXY}, $\Phi(f)$ is an element of $\MQXY$, and due to $X = \cup f[ Z ]$, the set $X$ does not contain a maximal point of $\Phi(f)$, hence $\Phi(f) \in \MstQXY$. On the other hand, for $R \in \MstQXY$, the set $Y$ contains all maximal points of $R$, and \eqref{inv_Phi} yields $X = \cup \; \inverse{\Phi(R)}[ Y ]$.

\EP

\begin{corollary} \label{coro_UXY_C}
Let $Z \subseteq Y$, $Q \in \fP(Z)$, and $R \in \UXY \cap \fC_Q(X \cup (Y \setminus Z),Z)$. If $Y \subseteq \downarrow_R Z$, then $R \in \fM_Q( X \cup (Y \setminus Z), Z )$.
\end{corollary}
\BP We have to show that $Z$ is an upper end of $R$. Suppose $z \in Z$ and $(z,x) \in R$ for an $x \in X \cup Y$. Then $x \inA \uparrow_R z \subseteqA \uparrow_R Z \subseteqA \uparrow_R Y$. But because $Y$ is an upper end of $R$, we have $ \uparrow_R Y = Y \subseteqA \downarrow_R Z$. There exists thus a $y \in Z$ with $(x,y) \in R$, and the convexity of $Q$ in $R$ yields $x \in Z$.

\EP

\begin{lemma} \label{lemma_tau}
The mapping
\begin{align} \nonumber
\tau_ {X,Y} : \UXY & \rightarrow \UXY, \\ \label{def_tau}
R & \mapsto \left( R \vert_X \right)^d \; \cup \; \left( ( X \times Y ) \setminus R \right) \; \cup \; \left( R \vert_Y \right)^d
\end{align}
is a well-defined self-inverse bijection. Moreover, for every subset $Z \subseteq Y$ and every p.o.r.\ $Q \in \fP(Z)$, the mapping $\tau_{X,Y}$ induces a self-inverse bijection between $\UXY \cap \fC_Q(X \cup (Y \setminus Z),Z)$ and $\UXY \cap \fC_{Q^d}(X \cup (Y \setminus Z),Z)$; in particular, for $Z = Y$, between $\MQXY$ and $\MQdXY$.
\end{lemma}
\BP Let $R \in \UXY$. The relation $T \equiv \tau_{X,Y}(R)$ is reflexive because it contains $( R \vert_X )^d$ and $( R \vert_Y )^d$. Due to $T \cap ( Y \times X ) = \emptyset$, we have
\begin{align} \label{impl_yX_xX}
y \in X & \; \Rightarrow \; x \in X \quad \mytext{for all }(x,y) \in T.
\end{align}
For $(x,y), (y,x) \in T$, we have thus $x \in X \Leftrightarrow y \in X$, and $x = y$ follows.

Let $(a,b), (b,c) \in T$. Implication \eqref{impl_yX_xX} leaves four possible cases:
\begin{itemize}
\item $c \in X$: Then $(a,b), (b,c) \in ( R \vert_X )^d$, hence $(a,c) \in T$.
\item $c \in Y$, $b \in X$: Then $(b,a) \in R \vert_X$ and $(b,c) \in ( X \times Y ) \setminus R$, hence $(b,c) \notin R$. The transitivity of $R$ yields $(a,c) \notin R$, thus $(a,c) \in T$.
\item $b, c \in Y$, $a \in X$: Then $(c,b) \in R \vert_Y$ and $(a,b) \in ( X \times Y ) \setminus R$, hence $(a,b) \notin R$. The transitivity of $R$ yields $(a,c) \notin R$, thus $(a,c) \in T$.
\item $a, b, c \in Y$: Then $(a,b), (b,c) \in ( R \vert_Y )^d$, hence $(a,c) \in T$.
\end{itemize}
Therefore, $T$ is a partial order relation on $X$, and due to $T \cap ( Y \times X ) = \emptyset$, the set $Y$ is an upper end in $T$. Because of $\tau_{X,Y}( \tau_{X,Y}(R)) = R$ for all $R \in \UXY$, the first part of the proposition is proven. 

Let $R \in \UXY \cap \fC_Q(X \cup (Y \setminus Z),Z)$  and $T \equiv \tau_{X.Y}(R)$, and let $a, b \in Z$, $y \in X \cup Y$, with $(a,y), (y,b) \in T$. In the case of $y \in X$, implication \eqref{impl_yX_xX} yields $a \in X$ in contradiction to $a \in Z \subseteq Y$. Therefore, $y \in Y$, hence $(y,a), (b,y) \in R \vert_Y$, and $y \in Z$ follows because $Z$ is convex in $R$. Therefore, $Z$ is convex in $T$, hence $T \in \UXY \cap \fC_{Q^d}(X \cup (Y \setminus Z),Z)$.

\EP

\begin{figure} 
\begin{center}
\includegraphics[trim = 75 630 210 70, clip]{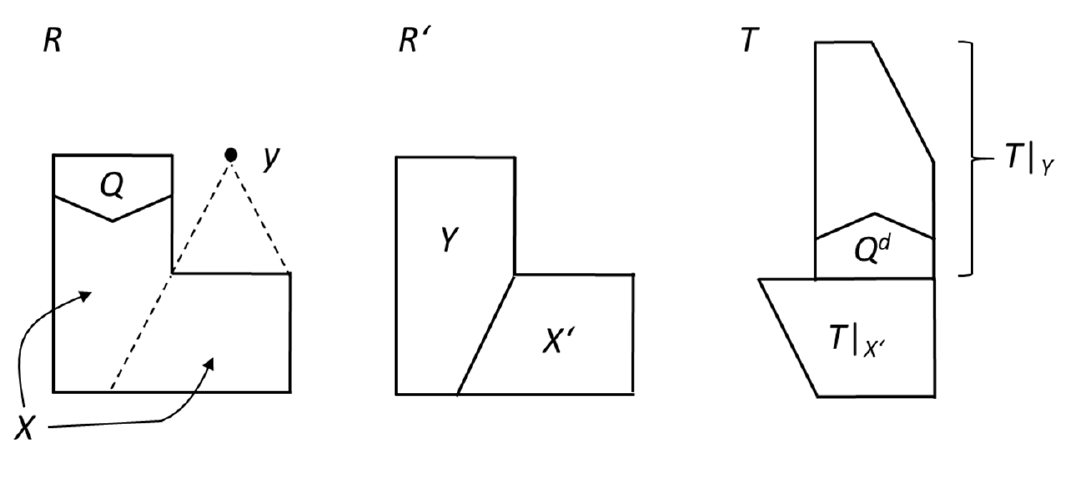}
\caption{\label{figure_sigmaR} The construction of $\sigma(R)$ in Theorem \ref{theo_Verne}. Explanations in text.}
\end{center}
\end{figure}

In the following theorem, we proof the first generalization of Ern\'{e}'s Theorem. The mapping $\sigma$ used in it is illustrated in Figure \ref{figure_sigmaR}; additional  examples are shown in Figure \ref{figure_Beispiel_Az}.

\begin{theorem} \label{theo_Verne}
Let $X$ and $Z$ be disjoint finite sets, $Q \in \fP(Z)$, let $y$ be a point not contained in $X \cup Z$, and let $A_y \equiv \{(y,y)\}$ be the antichain with carrier $\{y\}$. Then, with $W \equiv X \cup Z$,
\begin{align*}
\sigma : \MstQXZy & \rightarrow \CQdXZ \\
R & \mapsto \tau_{W \cap \downarrow_R y, \; W \setminus \downarrow_R y } ( R \vert_W )
\end{align*}
is a bijective mapping.
\end{theorem}
\BP Let $R \in \MstQXZy$, $R' \equiv R \vert_W$, $X' \equiv W \cap \downarrow_R y$, and $Y \equiv W \setminus \downarrow_R y$. The set $\Zy$ is an upper end in $R$ with $R \vert_\Zy = Q + A_y$. Therefore, $Z$ is an upper end in $R'$  with $R' \vert_Z = Q$. Furthermore, $Y$ is an upper end of $R'$, hence $R' \in \fU( X', Y )$, and 
\begin{align} \label{das_ist_T}
T & \equivA \tau_{X', Y} ( R' ) \eqA \left( R' \vert_{X'} \right)^d \; \cup \; \left( ( X' \times Y ) \setminus R' \right) \; \cup \; \left( R' \vert_Y \right)^d
\end{align}
is a well-defined element of $\fU( X', Y )$ according to Lemma \ref{lemma_tau}. $R \vert_\Zy = Q + A_y$ yields $Z \subseteq Y$, and because $Q + A_y$ is convex in $R$, $Q$ is convex in $R'$. Therefore, $Q^d = T \vert_Z$ is convex in $T$ according to the addendum in Lemma \ref{lemma_tau}.

We proceed with the proof of 
\begin{align} \label{Y_ZsigR}
Y & \eqA \uparrow_T Z.
\end{align}
``$\supseteq$'' is clear because $Y$ is an upper end in $T$ with $Z \subseteq Y$. Conversely, assume $x \in Y$. Because of $R \in \MstQXZy$, the set $\Zy$ contains all maximal points of $R$. There exists thus a maximal point $z \in \Zy$ of $R$ with $(x,z) \in R$. Due to $R \vert_\Zy = Q + A_y$, we must have $z \in Z \subseteq Y$, hence $(x,z) \in R' \vert_Y$, and $x \inA \uparrow_T Z$ follows.

Now we can show that $R$ is uniquely determined by $T$, i.e., that $\sigma$ is one-to-one. We have
\begin{align*}
\downarrow_R y & \eqA X' \cup \setyps \eqA ( W \setminus Y ) \cup \setyps \stackrel{\eqref{Y_ZsigR}}{\eqA} ( W \setminus \uparrow_T Z ) \cup \setyps, 
\end{align*}
hence $X' = W \setminus \uparrow_T Z$, and finally, Lemma \ref{lemma_tau} yields $R = \tau_{X',Y} ( \tau_{X',Y}( R ) ) =$ $\tau_{W \setminus \uparrow_T Z, \; \uparrow_T Z}( T )$.

It remains to show that $\sigma$ is onto. Let $R' \in \CQdXZ$ and $Y \equivA \uparrow_{R'} Z$, $X' \equiv X \setminus Y$. We have
\begin{align*}
R' & \; \inA \fU(X', Y ) \cap \fC_{Q^d}(X,Z) \\
& \eqA
\fU(X', Y ) \cap \fC_{Q^d}(X' \cup (Y \setminus Z),Z),
\end{align*}
and with $R \equiv \tau_{X', Y}(R')$, the addendum in Lemma \ref{lemma_tau} yields
\begin{align*}
R & \; \inA  \fU(X', Y ) \cap \fC_Q(X' \cup (Y \setminus Z),Z) \\
& \eqA \fU(X', Y ) \cap \fC_Q(X,Z).
\end{align*}

For every $x \in Y$, there exists a $z \in Z$ with $(z,x) \in R'$. Because of $Z \subseteq Y$, we have $(x,z) \in (R' \vert_Y)^d \subseteq R$. We conclude $Y \subseteq \; \downarrow_R Z$, and Corollary \ref{coro_UXY_C} delivers $R \in \fM_Q( X' \cup (Y \setminus Z), Z ) = \fM_Q( X, Z )$. With $R^o \equiv R \vert_X$, the set $( \downarrow_R z ) \setminus Z$ is a lower end of $R^o$ for every $z \in Z$. Moreover, because $Y$ is an upper end in $R$, the set $X'$ is a lower end in $R^o$. All together, the assignments
\begin{align*}
y & \mapsto X', \\
z & \mapsto ( \downarrow_R z ) \setminus Z \quad \mytext{for all } z \in Z,
\end{align*}
define a mapping $f \in \H_{Q + A_y}( \Zy, \L( R^o ) )$ with
\begin{align*}
\bigcup f \left[ \Zy \right] & \eqA \Big( ( \downarrow_R Z ) \setminus Z  \Big) \cup X' \\
& \; \supseteqA ( Y \setminus Z ) \cupA X' \eqA X,
\end{align*}
hence $f \in \FstQXZy$, because $\cup f \left[ \Zy \right] \subseteq X$ is trivial. With Corollary \ref{coro_MstQXY}, we conclude that the relation $\Phi(f)$ is an element of $\MstQXZy$. According to \eqref{def_Phi},
\begin{align*}
\Phi(f) & \eqA S(f) \; \cupA ( Q + A_y ) \; \cupA \bigcup_{x \in \Zy} ( f(x) \times \{ x \} ) \\
& \eqA R^o \; \cupA ( Q + A_y ) \; \cupA ( X' \times \setyps ) \; \cupA \bigcup_{z \in Z} ( ( ( \downarrow_R z ) \setminus Z ) \times \{ z \} ).
\end{align*}
Because of $\Phi(f) \vert_\Zy = Q + A_y$, the union on the right is a subset of $W \times W$, thus
\begin{align*}
\Phi(f) \vert_W & \eqA R^o \; \cupA Q \; \cupA \bigcup_{z \in Z} ( ( ( \downarrow_R z ) \setminus Z ) \times \{ z \} ).
\end{align*}
Due to $R^o = R \vert_X$, $Q = R \vert_Z$, and $\left( \left( \downarrow_R z \right) \cap Z \right) \times \{ z \} \subseteq Q$ for all $z \in Z$,
\begin{align*}
\Phi(f) \vert_W & \eqA R \vert_X \; \cupA R \vert_Z \; \cupA \bigcup_{z \in Z} ( ( \downarrow_R z ) \times \{ z \} ) \\
& \eqA R \vert_X \; \cupA R \vert_Z \; \cupA ( R \cap ( W \times Z ) ).
\end{align*}
Because of $R \in \fM_Q( X, Z )$, $Z$ is an upper end in $R$, therefore $R \cap (Z \times W ) \subseteq R \vert_Z$, and we arrive at $\Phi(f) \vert_W = R$. Thus, with $ \downarrow_{\Phi(f)} y = X' \cup \setyps $ and $W = X' \cup Y$,
\begin{align*}
\sigma( \Phi(f) ) & \eqA 
\tau_{W \cap \downarrow_{\Phi(f)} y, \; W \setminus \downarrow_{\Phi(f)} y } ( \Phi(f) \vert_W ) \\
& \eqA \tau_{X', Y} ( R ) \; \eqA \tau_{X', Y} ( \tau_{X', Y} ( R' ) ) \eqA R',
\end{align*}
because $\tau_{X', Y}$ is self-inverse according to Lemma \ref{lemma_tau}.

\EP

\section{An additional generalization step} \label{sec_VVerne}

Comparing Theorem \ref{theo_Verne} with our goal Theorem \ref{theo_VVErne_Intro}, we see that we have to skip the assumption about the convexity of $Q = R \vert_Z$. The required generalizations of the sets $\CQXY$ and $\MstQXY$ are the following sets:

\begin{definition} \label{def_NQXY}
Let $X$, $Y$ be disjoint sets and $Q \in \fP(Y)$. We define
\begin{align*}
\fI_Q( X, Y ) & \equivA \mysetdescr{ R \in \fP(X \cup Y) }{ R \vert_Y = Q }, \\
\fN^*_Q( X, Y ) & \equivA \mysetdescr{ R \in \fI_Q( X, Y ) }{ \max Q = \max R }.
\end{align*}
\end{definition}
Our leverage will be partial order relations $G$ with carrier $M \cup Y$, $M \subseteq X$, for which $G$ is the convex hull of $Y$, i.e., $\gamma_G(Y) = c(G) = M \cup Y$. With $\P(X)$ denoting the power-set of $X$, we define
\begin{definition} \label{def_GQAY}
Let $X$, $Y$ be disjoint sets and $Q \in \fP(Y)$. We define for all $M \in \P(X)$
\begin{align*}
\G_Q(M) & \equivA \mysetdescr{ G \in \fI_Q(M,Y)}{ \gamma_G(Y) = c(G) } \\
\mytext{and} \quad \fG_Q(X) & \equivA \bigcup_{M \in \P(X)} \G_Q(M).
\end{align*}
\end{definition}
We have $\G_Q(\emptyset) = \{ Q \}$ due to $\fI_Q(\emptyset, Y) = \{ Q \}$. If $Q$ is an antichain, then $\gamma_G(Y) = Y$ for  all $G \in \fI_Q(M,Y)$ and all $M \in \P(X)$, hence $\G_Q(M) = \emptyset$ for all $\emptyset \not= M \in \P(X)$.  If $Q$ is not an antichain, then there exist $a, b  \in Y$ with $(a,b) \in Q$, and for every $M \in \P(X)$, the transitive hull of $Q$ and $( \{ a \} \times M ) \cup ( M \times \{ b \} )$ is an element of $\G_Q(M)$.

If $Q$ is an antichain, then
\begin{align*}
\CQXY & \eqA \fI_Q( X, Y ), \\
\MstQXY & \eqA \fN^*_Q( X, Y ).
\end{align*}
The generalization of these equations is described in the following lemma:
\begin{lemma} \label{lemma_GQAY}
Let $X$, $Y$ be finite disjoint sets and $Q \in \fP(Y)$. The sets 
\begin{align*}
\fC_G( X \setminus c(G), c(G) ), & \quad G \in \fG_Q(X), \\
\fM^*_G( X \setminus c(G), c(G) ), & \quad G \in \fG_Q(X)
\end{align*}
form partitions of $\fI_Q(X,Y)$ and $\fN^*_Q(X,Y)$, respectively.
\end{lemma}
\BP Let $G \in \fG_Q(X)$. With $A$ being the antichain with carrier $X \setminus c(G)$, the p.o.r.\ $A \oplus G$ is an element of both, $\fC_G( X \setminus c(G), c(G) )$ and $\fM^*_G( X \setminus c(G), c(G) )$.

$R \in \fI_Q(X,Y)$ is clearly contained in $\fC_{R \vert_{\gamma_R(Y)}}( X \setminus \gamma_R(Y), \gamma_R(Y) )$. Let $R \in \fC_G( X \setminus c(G), c(G) ) \cap \fC_{G'}( X \setminus c(G'), c(G') )$ with $G, G' \in \fG_Q(X)$. There exist $M, M' \in \P(X)$ with $G \in \G_Q(M)$, and $G' \in \G_Q(M')$. Because of $G \subseteq R$ and the convexity of $c(G) = M \cup Y$ in $R$,
\begin{align*}
M \cup Y & \eqA \gamma_G(Y) \; \subseteqA \gamma_R(Y) \; \subseteqA M \cup Y,
\end{align*}
hence $\gamma_R(Y) = M \cup Y$. In the same way we see $\gamma_R(Y) = M' \cup Y$, thus $M = M'$, and $G = R \vert_{M \cup Y} = R \vert_{M' \cup Y} = G'$ follows.

Let $R \in \fN^*_Q(X,Y)$. Lemma \ref{lemma_ConstrConv} yields that $\gamma_R(Y)$ is an upper end of $R$ with
\begin{align*}
\max \gamma_R(Y) & \eqA \max Y \eqA \max Q \eqA \max R,
\end{align*}
and $R \in 
\fM^*_{R \vert_{\gamma_R(Y)}}( X \setminus \gamma_R(Y), \gamma_R(Y) )$ is shown. For $G, G' \in \fG_Q(X)$, we have
\begin{align*}
& R \inA \fM^*_G( X \setminus c(G), c(G) ) \; \capA \fM^*_{G'}( X \setminus c(G'), c(G') ) \\
\Rightarrow \quad & R \inA \fC_G( X \setminus c(G), c(G) ) \; \capA \fC_{G'}( X \setminus c(G'), c(G') ) \\
\Leftrightarrow \quad & G = G'.
\end{align*}

\EP

The following theorem completes the proof of Theorem \ref{theo_VVErne_Intro}:
\begin{theorem} \label{theo_VVErne}
Let $X$ and $Z$ be finite disjoint sets, $Q \in \fP(Z)$, and let $y$ be a point not contained in $X \cup Z$. Then
\begin{align*}
\# \fN^*_{Q + A_y}( X, \Zy ) & \eqA \# \fI_{Q^d}( X, Z ).
\end{align*}
\end{theorem}
\BP For every $G \in \fG_Q(X)$, Theorem \ref{theo_Verne} yields
\begin{align*}
\# \fM^*_{G + A_y}( X \setminus c(G), c(G) \cup \setyps ) & \eqA \# \fC_{G^d}( X \setminus c(G), c(G) ).
\end{align*}
Lemma \ref{lemma_GQAY} will deliver the result if we have shown
\begin{align*}
\fG_{Q^d}(X) & \eqA \mysetdescr{ G^d }{ G \in \fG_Q(X) }, \\
\fG_{Q + A_y}(X) & \eqA \mysetdescr{ G + A_y }{ G \in \fG_Q(X) }.
\end{align*}
The first equation is a direct consequence of the definition of $\G_Q(M)$, $M \in \P(X)$. For the second one, we get:

``$\subseteq$'': For $G' \in \fG_{Q + A_y}(X)$, there exists a $M \in \P(X)$ for which $G'$ is an element of $\fI_{Q + A_y}(M, \Zy)$, hence $G' \vert_{\Zy} = Q + A_y$ and $M \cup \Zy = \gamma_{G'}(\Zy)$. With $G \equiv G' \vert_{M \cup Z}$, we have $G \vert_Z = Q$. Let $x \in M \cup \Zy$. There exist $a, b \in \Zy$ with $(a,x), (x,b) \in G'$. $a = y = b$ yields $x = y$, $a, b \in Z$ yields $x \in \gamma_G(Z)$, and $a = y, b \in Z$ and $a \in Z, b = y$ both contradict $G' \vert_{\Zy} = Q + A_y$. We conclude $\gamma_G(Z) = \gamma_{G'}(\Zy) \setminus \setyps = M \cup Z$, thus $G \in \G_Q(M)$.

Let $m \in M$. Because of $\gamma_G(Z) = M \cup Z$, there exist $a, b \in Z$ with $(a,m), (m,b) \in G \subseteq G'$. $G' \vert_\Zy = Q + A_y$ yields $(m,y), (y,m) \notin G'$, hence
\begin{align*}
G' \cap ( M \times \setyps ) & \eqA \emptyset \eqA G' \cap ( \setyps \times M ),
\end{align*}
and we conclude $G' = G + A_y$.

``$\supseteq$'': For $G \in \fG_Q(X)$, we have $G \in \fI_Q(c(G) \setminus Z, Z)$ with $\gamma_G(Z) = c(G)$. But then $G + A_y \in \fI_{Q + A_y}(c(G) \setminus Z, \Zy)$ with $\gamma_{G + A_y}( \Zy ) = c(G) \cup \setyps = c(G + A_y)$, thus $G + A_y \in \fG_{Q + A_y}(X)$.

\EP

\begin{figure}
\begin{center}
\includegraphics[trim = 70 650 320 70, clip]{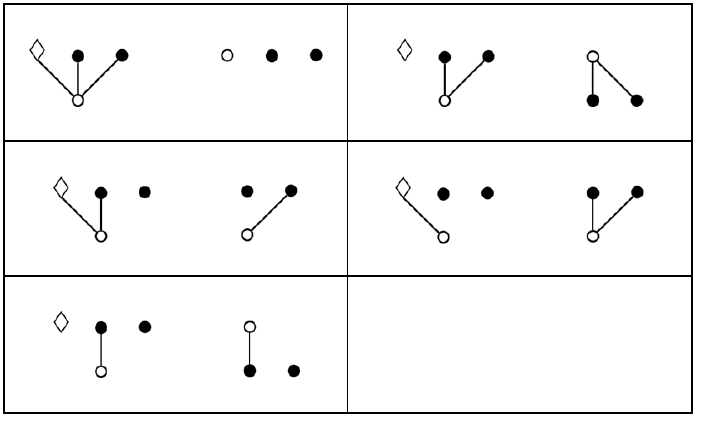}
\caption{\label{figure_Beispiel_Az} The Hasse diagrams of the partial order relations $R \in \fN^*_{A_2 + A_y}( X, \Zy ) = \fM^*_{A_2 + A_y}( X, \Zy )$ and of the corresponding partial order relations $\sigma(R) \in \fI_{A_2}( X \setminus \gamma_R(Z), \gamma_R(Z) ) = \fC_{A_2}( X, Z )$ for $\# X = 1$. For two of the pairs, two isomorphic copies exist.}
\end{center}
\end{figure}

\begin{figure}
\begin{center}
\includegraphics[trim = 70 620 290 70, clip]{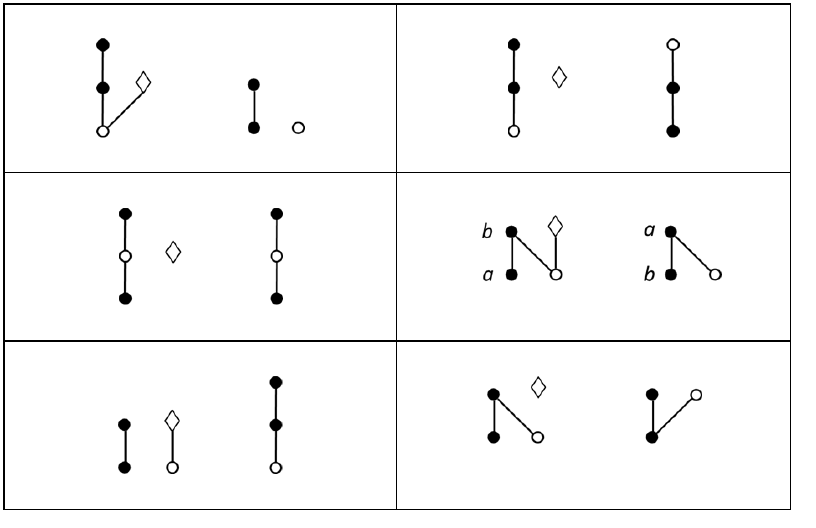}
\caption{\label{figure_Beispiel_Cz} The Hasse diagrams of the partial order relations $R \in \fN^*_{C_2 + A_y}( X, \Zy )$ and of the corresponding partial order relations $\sigma(R) \in \fI_{C_2^d}( X \setminus \gamma_R(Z), \gamma_R(Z) )$ for $\# X = 1$.}
\end{center}
\end{figure}

\begin{figure}
\begin{center}
\includegraphics[trim = 70 530 290 70, clip]{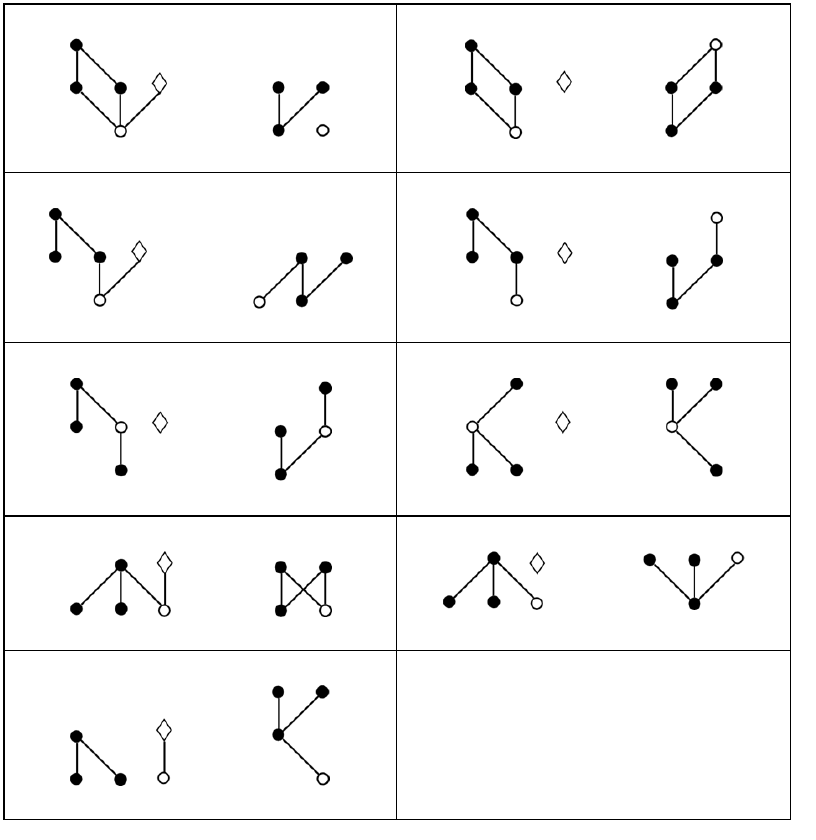}
\caption{\label{figure_Beispiel_V} The Hasse diagrams of the partial order relations $R \in \fN^*_{\Lambda + A_y}( X, \Zy )$ and of the corresponding partial order relations $\sigma(R) \in \fI_{\Lambda^d}( X \setminus \gamma_R(Z), \gamma_R(Z) )$ for $\# X = 1$. For three of the pairs, two isomorphic copies exist.}
\end{center}
\end{figure}

The Figures \ref{figure_Beispiel_Az}, \ref{figure_Beispiel_Cz}, and \ref{figure_Beispiel_V} show examples for the correspondence between the sets $\fN^*_{Q + A_y}( X, \Zy )$ and $\fI_{Q^d}( X, Z )$ for $\# X = 1$ and $Q \in \{ A_2, C_2, \Lambda \}$; here, $A_2$ is the two-element antichain, $C_2$ is the two-element-chain, and $\Lambda$ is the p.o.r.\ with $\Lambda$-shaped Hasse-diagram. In all three figures, the elements of $Q$ are shown as black dots, the element of $X$ is represented by a small circle, and the hollow diamond is the symbol for $y$.

\clearpage
\begin{thebibliography}{xx}
{\small
\bibitem{Alexandroff_1935} P. Alexandroff: Sur les espaces discrets. {\em Comptes Rendus de l'Académie des Sciences} {\bf 200} (1935), 1469--1471.     
      
\bibitem{Alexandroff_1937} P. Alexandroff: Diskrete R\"aume.
{\em Mat. Sb.} (N.S.) {\bf 2} (1937), 501--518.
      
\bibitem{Bayoumi_etal_1994} B. I. Bayoumi, M. H. El-Zahar, and S. M. Khamis: Counting two-dimensional posets. {\em Discrete Math.}  {\bf 131} (1994), 29--37.

\bibitem{Benoumhani_2006} M. Benoumhani: The number of topologies on a finite set. {\em J. Integer Seq.}  {\bf 9} (2006), Article 06.2.6.

\bibitem{Benoumhani_Kolli_2010} M. Benoumhani and M. Kolli: Finite topologies and partitions. {\em J. Integer Seq.}  {\bf 13} (2010), Article 10.3.5.

\bibitem{Birkhoff_1967} G. Birkhoff, {\em Lattice Theory}, Proc. Amer. Math. Soc. Coll. Publ. {\bf 25}, 3\textsuperscript{rd} ed.\ (1967).

\bibitem{Brinkmann_McKay_2002} G. Brinkmann and B. D. McKay: Posets on up to 16 points. {\em Order}  {\bf 19} (2002), 147--179.

\bibitem{Brinkmann_McKay_2005} G. Brinkmann and B. D. McKay: Counting unlabelled topologies and transitive relations. {\em J. Integer Seq.} {\bf 8} (2005), Article 05.2.1.

\bibitem{Butler_1972} K. K.-H. Butler: The number of partially ordered sets. {\em J. Combin. Theory (B)}  {\bf 13} (1972), 276--289.

\bibitem{Butler_Markowsky_1973} K. K.-H. Butler and G. Markowsky: Enumeration of finite topologies. in: {\em Proc. 4-th. Southeastern Conf. on Combinatorics, Graph Theory and Computing (Winnipeg)} (1973), 169–-184.

\bibitem{aCampo_2018} F. a Campo: A framework for the systematic determination of the posets on $n$ points with at least $\tau \cdot 2^n$ downsets. {\em Order} {\bf 36} (2019), 119--157. Published Online May 29, 2018, https://doi.org/10.1007/s11083-018-9459-2.

\bibitem{aCampo_Erne_2019} F. a Campo and M. Ern\'{e}: Exponential functions of finite posets and the number of extensions with a fixed set of minimal points. {\em J.\ Comb.\ Math.\ and Comb.\ Calc.\ } {\bf 110} (2019), 125--156.

\bibitem{Culberson_Rawlins_1991} J. C. Culberson and J. E. Rawlins: New results from an algorithm for counting posets.
{\em Order} {\bf 7} (1991), 361--374. 

\bibitem{Dershowitz_Zaks_1980} N. Dershowitz and S. Zaks: Enumerations of ordered trees. {\em Discrete Math.}  {\bf 31} (1980), 9--28.

\bibitem{Erne_1972} M. Ern\'{e}: {\em Struktur- und Anzahlformeln f\"{u}r Topologien auf endlichen Mengen.} PhD Dissertation, Westf\"{a}lische Wilhelms-Universit\"{a}t zu M\"{u}nster (1972).

\bibitem{Erne_1974} M. Ern\'{e}: Struktur- und Anzahlformeln f\"{u}r Topologien auf endlichen Mengen. {\em Manuscripta Math.}  {\bf 11} (1974), 221-–259.

\bibitem{Erne_1981} M. Ern\'{e}: On the cardinalities of finite topologies and the number of antichains in partially ordered sets. {\em Discrete Math.}  {\bf 35} (1981), 119--133.

\bibitem{Erne_1992} M. Ern\'{e}: The number of partially ordered sets with more points than incomparable pairs. {\em Discrete Math.}  {\bf 105} (1992), 49--60.

\bibitem{Erne_etal_2002} M. Ern\'{e}, J. Heitzig, and J. Reinhold: On the number of distributive lattices. {\em Electronic J. Comb.} {\bf 9} (2002), \#R24.

\bibitem{Erne_Stege_1991} M. Ern\'{e} and K. Stege: Counting finite posets and topologies. {\em Order}  {\bf 8} (1991), 247--265.

\bibitem{Erne_Stege_1995} M. Ern\'{e} and K. Stege: Combinatorial applications of ordinal sum decompositions. {\em Ars Combinatoria}  {\bf 40} (1995), 65--88.

\bibitem{Heitzig_Reinhold_2000} J. Heitzig and J. Reinhold: The number of unlabeled orders on fourteen elements. {\em Order}  {\bf 17} (2000), 333--341.

\bibitem{Heitzig_Reinhold_2002} J. Heitzig and J. Reinhold: Counting finite lattices. {\em Algebra Universalis}  {\bf 48} (2002), 43--53.

\bibitem{Kolli_2007} M. Kolli: Direct and elementary approach to enumerate topologies on a finite set. {\em J. Integer Seq.}  {\bf 10} (2007), Article 07.3.1.

\bibitem{Kolli_2014} M. Kolli: On the cardinality of the $T_0$-topologies on a finite set. {\em Intern. J. of Combin.}  {\bf } (2014), Article ID 798074. http://dx.doi.org/10.1155/2014/ 798074.

\bibitem{Parchmann_1975} R. Parchmann: On the cardinalities of finite topologies. {\em Discrete Math.}  {\bf 11} (1975), 161--172.

\bibitem{Pfeiffer_2004} G. Pfeiffer: Counting transitive relations. {\em J. Integer Seq.}  {\bf 7} (2004), Article 04.3.2.

\bibitem{Ragnarsson_Tenner_2010} K. Ragnarsson and B. E. Tenner: Obtainable sizes of topologies on finite sets. {\em J. Combin. Theory (A)} {\bf 117} (2010), 138--151.

\bibitem{Sharp_1966} H. Sharp: Quasi-orderings and topologies on finite sets. {\em Proc. Amer. Math. Soc.} {\bf 17} (1966), 1344--1349.

\bibitem{Sharp_1968} H. Sharp: Cardinality of finite topologies. {\em J. Combin. Theory}  {\bf 5} (1968), 82--86.

\bibitem{Stanley_1971} R. P. Stanley: On the number of open sets of finite topologies. {\em J. Combin. Theory (A)}  {\bf 10} (1971), 74-79.

\bibitem{Stanley_1974} R. P. Stanley: Enumeration of posets generated by disjoint unions and ordinal sums. {\em Proc. Amer. Math. Soc.}  {\bf 45-2} (1974), 295--299.

\bibitem{Vollert_1987} U. Vollert: {\em M\"{a}chtigkeiten von Topologien auf endlichen Mengen und Cliquenzahlen in endlichen Graphen.} PhD Thesis, University of Hannover, 1987.  
}
\end{thebibliography}
\end{document}